\documentclass[letterpaper, 12 pt]{article}
\usepackage{amsmath, amssymb, amsfonts, doublespace, theorem, epsfig}

\newcommand{\setleftmargin}[1]{
        \addtolength{\textwidth}{\oddsidemargin}
        \addtolength{\textwidth}{1in}
        \addtolength{\textwidth}{-#1}
        \setlength{\oddsidemargin}{-1in}
        \addtolength{\oddsidemargin}{#1}
        \setlength{\evensidemargin}{\oddsidemargin}
}
\newcommand{\setrightmargin}[1]{
        \setlength{\textwidth}{8.5in}
        \addtolength{\textwidth}{-\oddsidemargin}
        \addtolength{\textwidth}{-1in}
        \addtolength{\textwidth}{-#1}
}
\setleftmargin{1 in}
\setrightmargin{1 in}

\newcommand{\Z}{\mathbb{Z}}

\newcommand{\R}{\mathbb{R}}
\newcommand{\C}{\mathbb{C}}

\newcommand{\Cn}{\C^n}

\DeclareMathOperator{\Spec}{Spec}
\DeclareMathOperator{\Proj}{Proj}

\newcommand{\cs}{\C^\times}

\renewcommand{\L}{\mathcal{L}}
\renewcommand{\mod}{/\!\!/}
\newcommand{\xss}{X^{ss}}
\newcommand{\la}{\lambda}
\newcommand{\g}{\mathfrak{g}}
\newcommand{\tn}{\mathfrak{t}^n}
\renewcommand{\t}{\mathfrak{t}}
\renewcommand{\a}{\alpha}
\newcommand{\td}{\t^*}
\newcommand{\tdz}{\t^*_{\Z}}
\newcommand{\tdr}{\t^*_{\R}}
\newcommand{\D}{\Delta}

\newtheorem{theorem}{Theorem}[section]

{\theorembodyfont{\rmfamily}
\theoremstyle{plain}
\newtheorem{definition}[theorem]{Definition}
\newtheorem{example}[theorem]{Example}

\newtheorem{remark}[theorem]{Remark}
}

\newcommand{\qed}{\hfill \mbox{$\Box$}\medskip\newline}
\newenvironment{proof}{\noindent {\bf Proof:}}{\qed \par}

\begin{document}
\begin{spacing}{1.1}

\noindent
{\Large \bf Geometric invariant theory and projective toric varieties}
\bigskip\\
{\bf Nicholas Proudfoot}\footnote{Supported
by the Clay Mathematics Institute Liftoff Program
and the National Science Foundation Mathematical Sciences Postdoctoral Research Fellowship.}\\
Department of Mathematics, University of Texas,
Austin, TX 78712
\bigskip
{\small
\begin{quote}
\noindent {\em Abstract.}
We define projective GIT quotients, and 
introduce toric varieties from this perspective.
We illustrate the definitions by exploring the relationship
between toric varieties and polyhedra.
\end{quote}
}
\bigskip

\noindent
Geometric invariant theory (GIT) is a theory of quotients
in the category of algebraic varieties.
Let $X$ be a projective variety with ample line bundle $\L$,
and $G$ an algebraic group acting on $X$, along with a lift
of the action to $\L$.  The
GIT quotient of $X$ by $G$ is again a projective variety,
along with a given choice of ample line bundle.
With no extra work, we can consider varieties which are
projective over affine,
that is varieties can be written in the form $\Proj R$
for a reasonable graded ring $R$.
The purpose of this note is to give two equivalent definitions
of projective GIT quotients, one algebraic in terms of the homogeneous
coordinate ring $R$, and one more geometric,
and to illustrate these definitions with toric varieties.

A toric variety may be defined abstractly to be a normal variety
that admits a torus action with a dense orbit.
One way to construct such a variety is to take a GIT quotient
of affine space by a linear torus action, and it turns
out that every toric variety which is projective
over affine arises in this manner.
Given the data of a torus action 
on $\C^n$ along with a lift to the trivial line bundle, we define a polyhedron,
which will be bounded (a polytope) if and only if the corresponding toric
variety is projective.
We then use this polyhedron
to give two combinatorial descriptions
of the toric variety, one in the language of algebra
and the other in the language of geometry.

Much has been written about toric varieties, from many different
perspectives.  The standard text on the subject by Fulton \cite{Fu}
focuses on the relationship between toric varieties and fans.
The main difference between this approach and the one that we adopt here
is that a fan corresponds to an abstract toric variety, while
a polyhedron corresponds to a toric variety along with a choice of ample
line bundle.  In particular, there exist toric varieties which are
not projective over affine, and which are therefore do not come from 
polyhedra.  Since the primary purpose of this note is to introduce
projective GIT quotients, we will avoid fans altogether.
For an account of the development of the GIT perspective, including
a detailed explanation of its relationship to the fan construction,
one may consult the excellent survey paper by Cox \cite{Co}.
Toric varieties are used to illustrate a categorical perspective on
invariant theory in Dolgachev's book \cite[\S 12]{Do}, and they are studied
from the vantage point of multigraded commutative algebra in \cite[\S 10]{MS}.
Each of these three references covers most of the material that
is included in this paper and a lot more;
what we lack in depth and generality we hope to make up for with 
brevity and concreteness.

\paragraph{\bf Acknowledgments.}  I would like to express my gratitude
to Herb Clemens, Rob Lazarsfeld, and Ravi Vakil for organizing the 
conference at Snowbird out of which this paper grew.  
Also to the referee and David
Cox, both of whom provided valuable feedback on various drafts.

\begin{section}{Geometric invariant theory}\label{git}
Consider a graded noetherian algebra
$$R = \bigoplus_{m=0}^{\infty}R_m$$ which is finitely generated as an algebra over $\C$.
% The grading on $R$ induces an action of the multiplicative group $\cs$ on $\Spec R$
% by the formula $\lambda\cdot f = \lambda^mf$ for all $f\in R_m$.  The
% $\cs$-equivariant projection $R\onto R_0$ induces the inclusion $\Spec R_0\hookto\Spec R$
% of the $\cs$ fixed point set into $\Spec R$.  We then define
% $$\Proj R := \big(\Spec R\smallsetminus\Spec R_0\big)\big/\cs.$$
The variety $X = \Proj R$ is projective over the affine variety $\Spec R_0$, and comes
equipped with an ample line bundle $\L = \mathcal{O}_X(1)$.
% This variety is projective over $\Spec R_0$, and comes equipped with an ample line bundle
% \begin{eqnarray*}
% \L = \mathcal{O}_X(1) &:=& \big(\Spec R\smallsetminus\Spec R_0\big)\times_{\cs}\C\\
% &=& \Big(\big(\Spec R\smallsetminus\Spec R_0\big)\times\C\Big)\Big/\cs,
% \end{eqnarray*}
% where $\cs$ acts on $\C$ with weight $-1$.
Furthermore, by \cite[Ex. 5.14(a)]{Ha}, the integral closure (or normalization) of $R$
is isomorphic to the ring 
$$R' = \bigoplus_{m=0}^{\infty}\Gamma(X,\L^{\otimes m}).$$
The most important example for our applications will be the following.

\begin{example}\label{affine}
Let $R = \C[x_1,\ldots,x_n,y_0,\ldots,y_k]$, with $\deg x_i = 0$ and $\deg y_j = 1$.
Then $\Proj R \cong \C^n\times\C P^k$, and $\L$ restricts to the antitautological bundle
$\mathcal O(1)$ on $\{z\}\times\C P^k$ for all $z\in\C^n$.  
\end{example}

Let $G$ be a reductive algebraic group.  A good reference for general reductive groups is \cite{FH},
however the only groups that we will need for our applications in Section \ref{toric} are
subgroups of the algebraic torus $(\cs)^n$.
Suppose that we are given an action of $G$ on $R$ that preserves the grading; such an action induces
% This data is equivalent to an action of $G$ on $\Spec R$ which commutes with the action of $\cs$, or 
an action of $G$ on $X = \Proj R$ along with a lift of this action to the line bundle $\L$.
This lift is sometimes referred to as a {\em linearization} of the action of $G$ on $X$.

\begin{definition}\label{def}
The variety $X\mod G := \Proj(R^G)$ is called the GIT quotient of $X$ by $G$,
where $R^G$ is the subring of $R$ fixed by $G$.
\end{definition}

The reader is warned that, while the action of $G$ on $\L$ is not incorporated into the notation $X\mod G$,
it is an essential part of the data.  In particular, we will see in Example \ref{lins} that
it is possible to change the linearization of an action and obtain a vastly different quotient.

Definition \ref{def} is very easy to state, but not so easy to understand geometrically.
Our next goal will be to give a description of $X\mod G$ which depends more transparently
on the structure of the $G$-orbits in $X$.
Let $\L^*$ denote the line bundle dual to $\L$.

\begin{definition}\label{ss}
A closed point $x\in X$ is called {\em semistable} if, for all nonzero covectors $\ell\in\L^*_x$ over $x$,
the closure of the $G$-orbit $G\cdot (x,\ell)$ in the total space of $\L^*$ is disjoint from the zero section.
A point which is not semistable is called {\em unstable}.
The locus of semistable points will be denoted $\xss$.
\end{definition}

\begin{theorem}\label{equiv}
There is a surjective map $\pi:\xss\to X\mod G$, with $\pi(x) = \pi(y)$
if and only if the closures of the orbits $G\cdot x$ and $G\cdot y$
intersect in $\xss$.
\end{theorem}

\begin{proof}
We will provide only a sketch of the proof of Theorem \ref{equiv}; 
for a more thorough argument, see \cite[Prop 8.1]{Do}.
The projection from $R$ onto $R_0$ induces an inclusion of $\Spec R_0$ into $\Spec R$,
and the complement $\Spec R \smallsetminus \Spec R_0$ fibers over $\Proj R$ with fiber $\cs$.
The inclusion of $R^G$ into $R$ induces a surjection $\tilde\pi:\Spec R\to\Spec R^G$.
Let $x$ be an element of $X$, and let $\tilde x$ be a lift of $x$ to $\Spec R\smallsetminus\Spec R_0$.
Then
\begin{eqnarray*}
x\in\xss &\iff& \exists\hspace{3pt}\text{ a $G$-invariant section of $\L^{\otimes m}$
not vanishing at $x$ for some $m>0$}\\
&\iff& \exists \hspace{6pt}f\in R^G_m\text{ not vanishing at $\tilde x$ for some $m>0$}\\
&\iff& \tilde\pi(\tilde x)\notin\Spec R^G_0\subseteq\Spec R^G,
\end{eqnarray*}
hence $\tilde\pi(\tilde x)$ descends to an element of $\Proj R^G$.
Since the inclusion of $R^G$ into $R$ respects the gradings on the two rings,
this element does not depend on the choice of lift $\tilde x$,
hence $\tilde\pi$ induces a surjection $\pi:\xss\to\Proj R^G = X\mod G$.
Two points $x,y\in\xss$ with lifts $\tilde x,\tilde y\in\Spec R$
lie in different fibers of $\pi$ if and only there exists a $G$-invariant
function $f\in R^G_{>0}$ that vanishes at $\tilde x$ but not at $\tilde y$,
which is the case if and only if the closures of the $G$ orbits through $x$ and $y$
in $\xss$ are disjoint.
\end{proof}

Our proof of Theorem \ref{equiv} suggests that the variety $\Proj R$ may {\em itself} be
interpreted as a GIT quotient of $\Spec R$ by the group $\cs$.
Indeed, the 
grading on $R$ defines an action of $\cs$ on $R$
by the formula $\lambda\cdot f = \lambda^mf$ for all $f\in R_m$, and this
induces an action of $\cs$ on $\Spec R$.
Consider the lift of this action to the trivial line bundle $\Spec R \times \C$ given by letting $\cs$
act on the second factor by scalar multiplication.  The unstable locus for this linearized action
is exactly the subvariety $\Spec R_0\subseteq\Spec R$, and $\Proj R$ is the quotient
of $\Spec R\smallsetminus\Spec R_0$ by $\cs$.  This provides a geometric explanation of the irrelevance
of the irrelevant ideal in the standard algebraic definition of $\Proj$.

We conclude the section with an example that illustrates the dependence
of a GIT quotient on the choice of linearization of the $G$ action on $X$.

\begin{example}\label{lins}
As in Example \ref{affine}, let $R = \C[x_1,\ldots,x_n,t]$,
with $\deg x_i = 0$ for all $i$ and $\deg t = 1$.
Then $X\cong \C^n$, and $\L$ is trivial.
Let $G = \cs$ act on $R$ by the equations
$$\la\cdot x_i = \la x_i\hspace{8pt}\text{and}\hspace{8pt}
\la\cdot t = \la^{\alpha} t$$
for some $\alpha\in\Z$.  Geometrically, $G$ acts by scalar
multiplication, and $\alpha$ defines the linearization.
This action is not to be confused with the action of $\cs$ 
on $R$ given by the grading.

\vspace{\baselineskip}
\noindent{\em Case 1: $\alpha \geq 1$.}
In this case, $R^G = \C$, and $X\mod G = \Proj R^G$ is empty.
For every element $(x,\ell)\in\L^*$, we have 
$\displaystyle\lim_{\la\to 0}\la\cdot(x,\ell) = (0,0)$, hence
every $x\in X$ is unstable.

\vspace{\baselineskip}
\noindent{\em Case 2: $\alpha = 0$.}
With the trivial linearization of the $G$ action on $X$,
we have $R^G = \C[t]$, hence $X\mod G = \Proj R^G$ is a point.
The $G$ orbits in $\L^*$ are all horizontal, hence {\em every} point is semistable.
Since every $G$ orbit in $X$ contains the origin of $\Cn$ in its closure, Theorem \ref{equiv}
confirms that the quotient is a single point.

\vspace{\baselineskip}
\noindent{\em Case 3: $\alpha = -1$.}
In this case, $R^G = \C[x_1t,\ldots,x_nt]$ is a polynomial ring generated
in degree $1$, hence $X\mod G = \Proj R^G \cong \C P^{n-1}$.
We have $\la\cdot(x,\ell) = (\la x,\la^{-1}\ell)$, which limits to an element
of the zero section of $\L^*$ if and only if $x=0$.  Thus 
$\xss = \Cn\smallsetminus\{0\}$, and all $G$ orbits in $\xss$ are closed,
hence the GIT quotient is isomorphic to the quotient of $\xss$ by $G$ in the ordinary topological sense.

\vspace{\baselineskip}
\noindent{\em Case 4: $\alpha < -1$.}
In this case we still get $X\mod G \cong \C P^{n-1}$, but we now
obtain $\C P^{n-1}$ in its $(-\alpha)$-uple Veronese embedding.
\end{example}

Note that in Example \ref{lins},
multiplying $\alpha$ by a positive integer $m$ corresponds
to replacing the $G$-equivariant line bundle $\L$ on $X = \Cn$ with its
$m^{\text{th}}$ tensor power.
In general, this operation will have the effect of replacing the resulting ample line bundle 
on the GIT quotient $X\mod G$ by {\em its}
$m^{\text{th}}$ tensor power, as well (as we saw in Case 4).
\end{section}

\begin{section}{Toric varieties}\label{toric}
In this section we introduce and analyze toric varieties, which we will
think of as generalizations of Example \ref{lins} to higher dimensional tori.
As in the previous section, we let $$X = \Cn = \Proj\C[x_1,\ldots,x_n,t],$$ with $\deg x_i=0$
and $\deg t = 1$.  Fix an $n$-tuple $\alpha = (\alpha_1,\ldots,\alpha_n)$ of integers,
and let $T^n = (\cs)^n$ act on $R = \C[x_1,\ldots,x_n,t]$ by the equations
$$\la\cdot x_i = \la_ix_i\hspace{8pt}\text{and}\hspace{8pt}
\la\cdot t = \la_1^{\alpha_1}\ldots\la_n^{\alpha_n}\hspace{2pt}t$$
for $\la = (\la_1,\ldots,\la_n)\in T^n$.
Thus we have the standard coordinate action of $T^n$ on $\Cn$,
with a linearization to the trivial bundle given by $\alpha$.

\begin{definition}\label{toricdef}
A toric variety is a GIT quotient of $X$ by an algebraic subgroup 
$G\subseteq T^n$ for some $n$.
\end{definition}

A toric variety $X\mod G$ admits an action of the torus $T = T^n/G$ with
a single dense orbit.  A more standard approach to toric geometry is to define a toric
variety to be a normal variety along with a torus that acts with a dense orbit,
and then to prove that every such variety which is 
projective over affine arises from the construction of 
Definition \ref{toricdef}.  For the strictly projective case,
see \cite[\S 3.4]{Fu}.

Consider the exact sequence
$$1\to G\to T^n\to T\to 1.$$
Differentiating at the identity, we obtain an exact sequence of complex 
Lie algebras
$$0\to\g\to\tn\to\t\to 0.$$
Let $\{e_1,\ldots,e_n\}$ be the coordinate vectors in $\tn$, and let $a_i$ be the image
of $e_i$ in $\t$.  The vector space $\t$ is equipped with an integer lattice 
$\t_{\Z} = \ker\big(\operatorname{exp}:\t\to T\big)$.  Its dual $\td$ therefore
inherits a dual lattice, as well as a canonical real part 
$\tdr = \tdz\otimes_{\Z}\R$.
We now define the polyhedron
$$\D = \left\{\hspace{2pt}p\in\tdr\mid p\cdot a_i\geq\a_i
\text{ for all }i\hspace{2pt}\right\},$$
a subset of the real vector space $\tdr$.

There is a deep and extensive interaction between the toric variety $X\mod G$ and the polyhedron $\D$.
The $T$ orbits on $X$, for example, are classified by the faces of $\D$, with faces of real dimension $i$
corresponding to orbits of complex dimension $i$.
If $\D$ is simple (exactly $\dim_\R \D$ facets meet at each vertex), then $X\mod G$ is an orbifold \cite{LT}, and the 
Betti numbers of $X\mod G$ are determined by the equation
\begin{equation}\label{betti}
\sum_{i=0}^{d} b_{2i}(X\mod G) \hspace{2pt}q^i = \sum_{i=0}^{d} f_i(\D) (q-1)^i,
\end{equation}
where $d = \dim_\C X\mod G = \dim_\R\D$, and $f_i(\D)$ is the number of faces of dimension $i$.  
This fact has been famously used by Stanley to
characterize the possible face vectors of simple polytopes \cite{St},
and can be proven in many ways.  
One beautiful (though unnecessarily technical) proof uses the Weil conjectures;
it amounts simply to observing that the right hand side of 
Equation \eqref{betti} may be interpreted
as the number of points on an $\mathbb{F}_q$ model of $X\mod G$.
For a more detailed discussion of Betti numbers, the Weil conjectures,
and Stanley's theorem, see \cite[\S 4.5 and 5.6]{Fu}.
In this note, we will content ourselves with using $\D$ to describe the invariant ring $R^G$
and the semistable locus $\xss$.

Let $C_\D$ be the cone over $\D$, that is
$$C_\D := \left\{(p,r)\in\tdr\times\R\mid r\geq 0\text{ and }p\in r\cdot\D\right\}^{cl},$$
where ${cl}$ denotes closure inside of $\tdr\times\R$.
The following figure illustrates the cone over an interval and the cone over a half line.
Note that the closure is necessary to include the positive $x$-axis in the cone over the half line.
\begin{figure}[h]
\centerline{\epsfig{figure=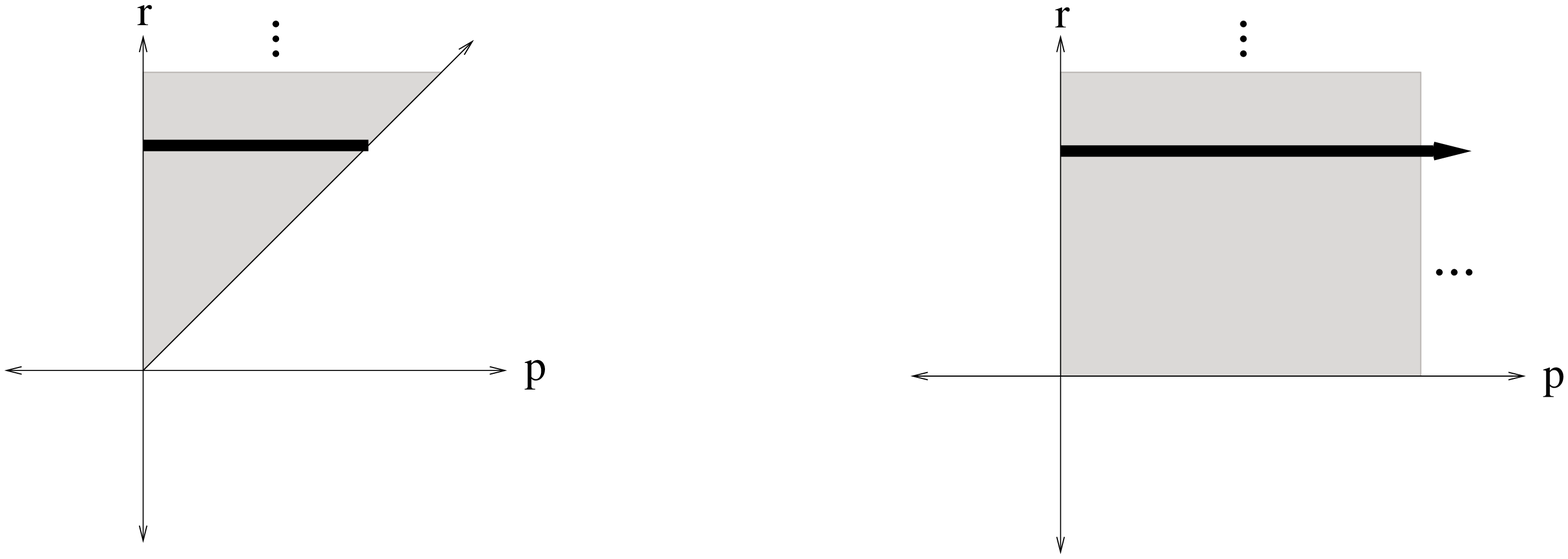, height=4cm}}
\end{figure}

Let $S_\D := C_\D\hspace{2pt}\cap\hspace{2pt}\big(\tdz\times\Z\big)$ be the semigroup consisting of all of the lattice points in $C$.
We may then define the semigroup ring $\C[S_\D]$, an algebra over $\C$
with additive basis indexed by the elements of $S_\D$, and multiplication
given by the semigroup law.  This ring has a non-negative integer
grading given on basis elements by the final coordinate of the corresponding
lattice points.
The following theorem provides a combinatorial interpretation of the homogeneous coordinate ring
$R^G$ of $X\mod G$.

\begin{theorem}\label{invt}
$R^G\cong \C[S_\D]$.
\end{theorem}

\begin{proof}
Suppose that we are given an element $(p,r)\in S_\D$,
and let $r_i = p\cdot a_i - r\a_i \in \Z_{\geq 0}$ for all $i$.
To this element, there corresponds a $G$-invariant monomial
$m_{(p,r)} = x_1^{r_1}\ldots x_n^{r_n}t^r\in R$.
This correspondence defines a bijection from $S_\D$ to the monomials of $R^G$,
which extends to a graded ring isomorphism $\C[S_\D]\cong R^G$.
\end{proof}

For all $i\in\{1,\ldots,n\}$, let $F_i = \{p\in\D\mid p\cdot a_i = \a_i\}$.
The set $F_i$ is the locus of points on $\D$ on which the 
$i^\text{th}$ defining linear form is minimized,
and therefore it is a face of $\D$. 
If $\a$ is chosen generically, then $F_i$ will either be a facet
or it will be empty.  In general, however, $F_i$ may be a face of any dimension.
The following theorem provides a combinatorial interpretation of the semistable locus $\xss$.

\begin{theorem}\label{xss}
For any point $x\in X$, let $A = \{i\mid x_i=0\}$ 
be the set of coordinates at which $x$ vanishes.
Then $x$ is semistable if and only if $\displaystyle\bigcap_{i\in A}F_i\neq\emptyset$.
\end{theorem}

\begin{proof}
In the proof of Theorem \ref{equiv}, we argued that $x$ is stable
if and only if there is a positive degree element of $R^G$ that does not vanish at $x$.
This will be the case if and only if there is $G$-invariant monomial
of positive degree which is
supported on the complement of $A$.  By Theorem \ref{invt}, $G$-invariant
monomials of degree $r$ correspond to lattice points in $r\cdot\D$,
and those that are supported on the complement of $A$ correspond to those
that lie on $r\cdot F_i$ for all $i\in A$.  Such a monomial exists
if and only if $\displaystyle\bigcap_{i\in A}F_i\neq\emptyset$.
\end{proof}

\begin{remark}
Given a polyhedron $\Delta\in\R^d$, there are infinitely many
different ways to present it as the set of solutions to a finite
set of affine linear inequalities.  Indeed, even the number of such
inequalities is not determined; it is only bounded below by the number
of facets of $\Delta$. 
Theorems \ref{invt} and \ref{xss} are valid regardless of the choice of
presentation of $\Delta$,
and the presentation is irrelevant to the application 
of Theorem \ref{invt}.  
On the other hand, if we want to apply Theorem \ref{xss} in an example,
it is essential to be given $\Delta$ along with its presentation.  We
can then set $n$ to be the number of defining inequalities,
and reconstruct the group $G\subseteq T^n$ by reversing the construction
that we gave for $\Delta$.
\end{remark}

\vspace{8pt}
We conclude the section by using Theorems \ref{invt} and \ref{xss} to
compute the toric varieties associated to an assortment of polytopes.
An integer vector $a\in\tdz$ is called {\em primitive} if it cannot be expressed
as a multiple of another element $a'\in\tdz$ by an integer greater than $1$.
In each of the following examples
we will implicitly assume that the given polytope is cut out by the minimum
possible number of linear forms $\{a_1,\ldots,a_n\}$ in the dual vector space,
and that each of these forms is a primitive integer vector.

\begin{example}\label{ray}
Let $\D = \R^+\subset\R$ be the set of non-negative real numbers.
Then $C_\D = (\R^+)^2$ and $\C[S_\D]\cong\C[x,t]$, with
$\deg x = 0$ and $\deg t = 1$.  This tells us that the associated
toric variety is $\Proj \C[x,t] = \C$.  Geometrically, we have $T^n = T = \cs$,
and $G$ is the trivial group, hence we are building
$\C$ as a trivial GIT quotient of $\C$ itself.
More generally, the toric variety associated to the positive orthant in $\R^d$
is $\C^d$, equipped with the trivial line bundle.
\end{example}

\begin{example}\label{interval}
Let $\D = [0,1]$ be the unit interval in $\R$.
Then $S_\D \cong\C[x,y]$ is a polynomial ring in two variables of degree $1$,
and the associated toric variety is $\C P^1$.  Geometrically,
$\cs$ acts by scalars on $\C^2$, and the origin is the unique unstable point.
More generally, the toric variety associated to the standard $d$-simplex 
in $\R^d$ is $\C P^d$ with its antitautological line bundle.
\end{example}

\begin{example}\label{biginterval}
Let $\D = [0,m]\subset\R$ for some positive integer $m$.
The action of $\cs$ on $\C^2$ and the semistable locus
are unchanged from Example \ref{interval}, hence the associated
toric variety is again $\C P^1$.  Its homogeneous coordinate ring
$\C[S_\D]$, however, is isomorphic to the subring of $\C[S_{[0,1]}]$
spanned by homogeneous polynomials in degrees which are multiples
of $m$.  Hence the line bundle that we obtain on $\C P^1$
is the $m^{\text{th}}$ tensor power of the antitautological line bundle.
More generally, dilating $\D$ by a positive integer $m$
corresponds to taking the $m^{\text{th}}$ tensor
power of the ample line bundle on the toric variety.
(See the end of Section \ref{git}.)
\end{example}

\begin{example}\label{square}
Let $\D = [0,1]\times [0,1]\subset \R^2$.
This corresponds to an action of $(\cs)^2$ on $\C^4$,
given in coordinates by $$(\la,\mu)\cdot(z_1,z_2,z_3,z_4)
=(\la z_1,\la z_2,\mu z_3,\mu z_4).$$
The unstable locus consists of the points where either
$z_1=z_2=0$ or $z_3=z_4=0$, and the quotient of the semistable
points by $(\cs)^2$ is isomorphic to $\C P^1\times\C P^1$.
On the algebraic side, we have 
$$\C[S_\D]\cong\C[x,y,z,w]/\langle xz-yw\rangle,$$
where $x,y,z,$ and $w$ are generators in degree $1$
corresponding to a cyclic ordering of the vertices of $\D$.
In general, the toric variety corresponding to the product
of two polytopes is isomorphic to the product of the corresponding
toric varieties, in the Segr\'e embedding.
\end{example}
\end{section}

\footnotesize{

}
\end{spacing}
\end{document}